\documentclass[11pt,a4paper] {amsart}
\newtheorem{thm}{Theorem}
\begin{document}
\title[A combinatorial formula for homogeneous moments]
{A combinatorial formula for \\
homogeneous moments}
\author{Michael G. Eastwood}
\address{School of Pure Mathematics, University of Adelaide, SA 5005,
Australia}
\email{meastwoo@maths.adelaide.edu.au}
\author{Nuno M. Rom\~ao}
\address{School of Pure Mathematics, University of Adelaide, SA 5005,
Australia}
\email{nromao@maths.adelaide.edu.au}
\thanks{The authors are supported by the Australian Research Council.}
\begin{abstract}We establish a combinatorial formula for homogeneous moments
and give some examples where it can be put to use.
An application to the statistical mechanics of interacting
gauged vortices is discussed.\end{abstract}
\maketitle

\section{Introduction}
In this article we shall prove the following.
\begin{thm} \label{thethm}
Suppose we are given $n$ integrable real-valued functions $J_1,J_2,\ldots,J_n$
on a measure space $M$, an integer $m\ge 0$, and a non-zero constant~$C$, so
that 
\begin{equation}\label{given}\int_M(v_1J_1+v_2J_2+\cdots+v_nJ_n)^{2m}=C,
\end{equation}
for all $(v_1,v_2,\ldots,v_n)\in{\mathbb{R}}^n$ with
$v_1{}^2+v_2{}^2+\cdots+v_n{}^2=1$. Then there is a rational number
\begin{equation}\label{formulae}
I_{m,n}=\frac{\Gamma(\frac12)\Gamma(m+\frac{n}2)}
{\Gamma(m+\frac12)\Gamma(\frac{n}2)}\end{equation}
(depending only on $m$ and~$n$) so that 
\begin{equation}\label{thesis}
\int_M(J_1{}^2+J_2{}^2+\cdots+J_n{}^2)^m=I_{m,n}\,C.
\end{equation}
\end{thm}
We refer to quantities such as the integrals in the left-hand side of (\ref{given}) and (\ref{thesis}) as {homogeneous moments} (in the given
functions $J_1, \ldots, J_n$).
For applications, we have in mind that $M$ should be a compact smooth manifold
equipped with a volume form and that $J_1,J_2,\ldots,J_n$ should be smooth
functions. A simple example may be given by taking $M=S^{n-1}$ with its usual
round metric and $J_1,J_2,\ldots,J_n$ to be the co\"ordinate functions for the
standard embedding $S^{n-1}\hookrightarrow{\mathbb{R}}^n$. In this case,
$v_1J_1+v_2J_2+\cdots+v_nJ_n:S^{n-1}\to{\mathbb{R}}$ is simply another
co\"ordinate function on the sphere and (\ref{given}) is evident by rotational
invariance. (The constant $C$ may always be realised by setting $v_n=1$ and all
other $v_j$'s equal to zero: $C=\int_MJ_n{}^{2m}$.) In this example, reflection
symmetry directly implies that the integral on the left-hand side of
(\ref{given}) would vanish if $2m$ were replaced by an odd integer; an easy
argument shows that this would also be true in the general case of our
hypothesis. Thus we are not losing generality by assuming that the degree of
homogeneity in (\ref{given}) is even.

Our motivation, however, stems from a more substantial example in which $M$ is
a K\"ahler manifold whose structure is invariant under the action
of~${\mathrm{SO}}(3)$. In this case $n=3$ and the three functions
$J_1,J_2,J_3:M\to{\mathbb{R}}$ are the components of the associated moment
map $M\to{\mathfrak{so}}(3)^*$. The action of ${\mathrm{SO}}(3)$ arises
because $M$ is the moduli space of gauged vortices on the round
$2$-sphere~\cite{ManSut}.
Further discussion of this case is provided at the end of this article. When
$n=3$, the rational numbers $I_{m,3}$ are, in fact, integers: $I_{m,3}=2m+1$.
That the $I_{m,n}$ are well-defined rational numbers is most easily proved
without the explicit formula~(\ref{formulae}). Having done this
in~\S\ref{Combinatorics}, we shall establish (\ref{formulae}) in
\S\ref{Formulae}. In~\S\ref{Examples}, we give a geometric setting for
our hypothesis (\ref{given}) with some examples, and conclude in \S\ref{Vortices} with a brief discussion of vortices on the $2$-sphere.\\

\section{Combinatorics}\label{Combinatorics} 
We may evidently extend (\ref{given}) by homogeneity to conclude that
\begin{equation}\label{polly}\int_M(v_1J_1+v_2J_2+\cdots+v_nJ_n)^{2m}=
C(v_1{}^2+v_2{}^2+\cdots+v_n{}^2)^m,\end{equation} 
for all $(v_1,v_2,\ldots,v_n)\in{\mathbb{R}}^n$. Also, by simultaneously
rescaling all the functions $J_1,J_2,\cdots,J_n$ we may suppose without loss of
generality that $C=1$ and our task is now to compute
$$I_{m,n}\equiv\int_M(J_1{}^2+J_2{}^2+\cdots+J_n{}^2)^m$$
from this information. To do this we may regard (\ref{polly}) as the equality
of two polynomials in the $u$-variables and deduce the equality of their
coefficients. Examining the left hand side of~(\ref{polly}), this means that we
may compute
$$\int_MJ_1{}^{r_1}J_2{}^{r_2}\cdots J_n{}^{r_n}\quad\mbox{for all }
r_1+r_2+\cdots+r_n=2m.$$
Indeed, since only integers are involved in expanding the two sides
of~(\ref{polly}), it follows that these integrals are all rational. This is
more than enough to compute $I_{m,n}$ and to conclude that these quantities are
also rational. 

Though the preceding argument is straightforward in principle, in practise it
is almost useless in establishing formulae for~$I_{m,n}$. We shall rectify
this deficiency in the following section. In the meantime, let us compute
$I_{3,3}$ by bare hands. We are supposing that 
$$\int_M(v_1J_1+v_2J_2+v_3J_3)^6=(v_1{}^2+v_2{}^2+v_3{}^2)^3$$
and, by expanding as polynomials in $v_1,v_2,v_3$, we conclude that 
$$\begin{array}c\int_MJ_1{}^6=\int_MJ_2{}^6=\int_MJ_3{}^6=1\\[8pt]
\int_MJ_1{}^4J_2{}^2=\int_MJ_1{}^2J_2{}^4=\cdots=
\int_MJ_2{}^2J_4{}^4=1/5\\[8pt]
\int_MJ_2{}^2J_2{}^2J_3{}^2=1/15\end{array}$$
(and more besides). Therefore, 
$$\begin{array}{rcl}
\int_M(J_1{}^2+J_2{}^2+J_3{}^2)^3&=&\int_M(J_1{}^6+\cdots
+3J_1{}^4J_2{}^2+\cdots+6J_1{}^2J_2{}^2J_3{}^2)\\[8pt]
&=&3\times 1+6\times 3\times 1/5+6\times 1/15=7.\end{array}$$
Such na\"{\i}ve computations of $I_{m,n}$ rapidly get out of hand for large $m$
and~$n$.\\

\section{Calculation of $I_{m,n}$} \label{Formulae}
Given that the argument we used above for the existence of $I_{m,n}$ is purely
combinatorial, it is clear that this quantity does not depend on the measure
space $M$. Thus we can calculate $I_{m,n}$ by evaluating both the integral and
the constant $C$ in (\ref{thesis}) for a specific model where the hypothesis
(\ref{given}) is satisfied. We shall do this for the first example mentioned in
the Introduction, where $M$ is $S^{n-1}$ with metric induced from the embedding
$\iota: S^{n-1} \hookrightarrow \mathbb{R}^n$ and 
$J_{j}:S^{n-1} \rightarrow [-1,1]$ are the standard cartesian co\"ordinates.

In this example, the integral on the left-hand side of (\ref{thesis}) is just
the well-known volume of $S^{n-1}$,
\begin{equation}\label{firstint}
\int_{S^{n-1}}({J_1}^2+\cdots+{J_{n}}^2)^{m}=\int_{S^{n-1}}1=
\frac{2 \pi^{n/2}}{\Gamma(\frac{n}{2})}.
\end{equation}
To evaluate the constant $C=\int_{S^{n-1}}{J_n}^{2m}$ on the right-hand side,
we start by using cylindrical co\"ordinates to write the euclidean metric $g$
on in $\mathbb{R}^{n}$ as
\begin{equation}\label{euclidean}
g={\rm d}x^{2}+{\rm d}r^2 + r^2 \,g_{S^{n-2}}.
\end{equation}
Here, $x=J_n$, $r=\sqrt{{J_1}^2+\cdots+{J_{n-1}}^2}$ and $g_{S^{n-2}}$ denotes
the metric on $S^{n-2}$. On $S^{n-1}\subset \mathbb{R}^n$, we have the relation
\[
x^2 + r^2 = 1 \quad \Rightarrow \quad x \,{\rm d}x = - r\,{\rm d}r,
\]
therefore (\ref{euclidean}) pulls back to $S^{n-1}$ as
\[
g_{S^{n-1}}=\iota^{*}g = \frac{{\rm d}x^2}{1-x^2}+(1-x^2)\, g_{S^{n-2}}.
\]
In these co\"ordinates, the volume form on $S^{n-1}$ can then be written as
\begin{equation}\label{generalhatbox}
{\rm dvol}_{S^{n-1}}=\sqrt{\det(g_{S^{n-1}})}\,{\rm d}x\wedge {\rm
dvol}_{S^{n-2}}
=(1-x^2)^{(n-3)/2}\,{\rm d}x \wedge {\rm dvol}_{S^{n-2}}.
\end{equation}
We observe in passing that the case $n=3$ is special, as (\ref{generalhatbox})
then gives an identification of the volume forms on $S^{2}$ and on the cylinder
in $\mathbb{R}^{3}$ defined by $r=1$ and $|x|\le 1$; this is the celebrated
hat-box theorem of Archimedes, of which our equation (\ref{generalhatbox}) may
be regarded as a generalisation to arbitrary dimensions. The constant on the
right-hand side of (\ref{thesis}) can now be evaluated as
\begin{eqnarray}
\int_{S^{n-1}}{J_{n}}^{2m}&=&\int_{[-1,1]\times S^{n-2}}
x^{2m}(1-x^2)^{(n-3)/2}\,{\rm d}x \wedge {\rm dvol}_{S^{n-2}} \nonumber \\
&=& \int_{S^{n-2}} 1 \times \int_{0}^{1}
t^{m-1/2} (1-t)^{(n-3)/2}\, {\rm d}t \nonumber \\
&=&\frac{2 \pi^{(n-1)/2}}{\Gamma(\frac{n-1}{2})}\, 
{\textstyle B\left(m+\frac{1}{2},\frac{n-1}{2}\right)} \nonumber \\
&=& \frac{2 \pi^{(n-1)/2} \, \Gamma(m+\frac{1}{2})}{\Gamma(m+\frac{n}{2})}, 
\label{secondint}
\end{eqnarray}
where $B$ denotes Euler's Beta-function~\cite{AbrSte}. Finally, we obtain
$I_{m,n}$ as a quotient of (\ref{firstint}) and (\ref{secondint}), which 
yields (\ref{formulae}) using $\Gamma(\frac{1}{2})=\sqrt{\pi}$.

We would like to make two remarks on the formula (\ref{formulae}):
\begin{enumerate}
\renewcommand{\labelenumi}{\Alph{enumi})}

\item Notice that (\ref{formulae}) gives $I_{0,n}=1$, and so it follows from
the recursion relation $\Gamma(z+1)=z\,\Gamma(z)$ that $I_{m,n}$ is obviously
a rational number for all $m\in \mathbb{N}$. More explicitly, we can write
\[
I_{m,n}=\left\{
\begin{array}{c@{\quad{\rm if}\;\,n\;\,{\rm is}\;}l}
{\displaystyle \frac{2^{2m-1}\left(m+\frac{n}{2}-1\right)!(m-1)!}
{(2m-1)! \left(\frac{n}{2}-1 \right)!} } & \mbox{even,}\\[3.5 ex]
{\displaystyle \frac{(2m+n-2)! (m-1)! \left( \frac{n-3}{2}\right)!}
{2(2m-1)! (n-2)! \left( m+ \frac{n-3}{2}\right)!} } & \mbox{odd.}
\end{array}
\right.
\]

\item Our result for $I_{m,n}$ as expressed in equation (\ref{formulae}) 
turns out to be a value of Gau\ss's hypergeometric function,
\[
I_{m,n}=
{\textstyle \left._2F_1 \left( 2m,n-1; m+\frac{n}{2};
\frac{1}{2}\right)\right.} ,
\]
which comes as a consequence of Gau\ss's second summation theorem
(cf.~\cite{SlapFq}, p.~32). The property $I_{m,n} \in \mathbb{Q}$ provides
examples of the curious fact that the hypergeometric function sometimes assumes
rational values when its argument and parameters are rational. This behaviour
is general for the geometric series, but it is not understood under which
circumstances it generalises to hypergeometric series.\\

\end{enumerate}

\section{Geometry and examples}\label{Examples}

We shall now describe a geometric setup where the hypothesis (\ref{given})
naturally arises, and which inspires generalisations of 
Theorem~\ref{thethm}. 
Let us suppose that $V$ is a real
inner product space and that a compact
Lie group $G$ has an orthogonal represention on $V$, 
$\rho: G \rightarrow {\rm O}(V)$.
Typically, we take as $M$ a smooth manifold where $G$ acts,
equipped with an invariant volume form, which may in turn be induced by
a Riemannian or symplectic structure on~$M$. 
We also assume that 
$J:M \rightarrow V$ is a $G$-equivariant mapping.
At this point, symmetry considerations will severely restrict the 
integrals
\begin{equation}\label{integrals}
\Phi_m(v):=\int_M\langle v,J\rangle^{2m}
\end{equation}
with $v \in V$. In fact, for all $g \in G$,
\[
\Phi_m (\rho(g)v)=\int_{M}\langle \rho(g)v,J\rangle^{2m}=
\int_{M}\langle v,\rho(g^{-1})J \rangle^{2m}=\Phi_{m}(v),
\] 
and this implies that $\Phi(v)$ is an invariant
polynomial restricted by Weyl's classical invariant theory~\cite{WeyCG}.
Because $G$ is represented by orthogonal transformations, $||v||^{2m}$
is one of the invariant polynomials of degree $2m$; however, these
ingredients are still not enough to enforce (\ref{given}).
In the following, we discuss
three realisations of this general setup that 
lead to condition (\ref{given}) being satisfied, or to more general
conditions that still allow us to use the arguments
in sections~\ref{Combinatorics} and~\ref{Formulae} to determine 
$\int_{M}\langle J , J\rangle^{2m}$ by combinatorial means. \\

\noindent
{\bf Example 1.}
Let $G={\rm SO}(V)$, $\rho$ be the defining representation, and
$J:M \hookrightarrow V$ be the inclusion of any invariant measurable set.
Then necessarily $\Phi_m(v)=C||v||^{2m}$. A particular case is when $M$ is
a single ${\rm SO}(V)$-orbit embedded in $V$, which yields the example 
that we used in the calculations of section~\ref{Formulae}.\\

\noindent
{\bf Example 2.}
Let $G={\rm SO}(4)$, $V=\mathfrak{so}(4)\cong \bigwedge^{2}\mathbb{R}^4$ 
and $\rho$ be the adjoint representation. The ring of invariants is 
freely generated by two polynomials of degree two~\cite{AudST}, 
the squared norm $||v||^2$ and the pfaffian ${\rm Pf}(v)$. 
Thus we can write
\[
\Phi_{m}(v)=\sum_{j=0}^{m}C_{j}\,||v||^{2(m-j)}\,{\rm Pf}(v)^{j}
\]
for suitable real constants $C_j$. Our hypothesis in the form (\ref{polly})
will hold at least for 
$v$ in the hypersurface  defined by ${\rm Pf}(v)=0$, and the argument in 
section~\ref{Combinatorics} will still lead to 
a combinatorial formula for the homogeneous moments in (\ref{thesis}). Whether
(\ref{polly}) holds more generally depends on the choice of $M$ and
$J$. For instance, one can show that, for $m=1$, (\ref{polly}) holds 
for all $v \in V$ if and only if, say,
\begin{equation} \label{orthogonality}
\int_{M} J_{12} J_{34} = 0,
\end{equation}
where the indices refer to the standard basis of $\bigwedge^{2}\mathbb{R}^{4}$. An interesting example is when $M$ is 
the 4-manifold
of simple 2-vectors in $\bigwedge^{2}\mathbb{R}^4$ of unit norm, the 
adjoint orbit of ${\rm SO}(4)$ given by the algebraic equations 
${\rm Pf}(v)=0$ and $||v||^2=1$, and 
$J:M\hookrightarrow \bigwedge^{2}\mathbb{R}^4$ is the inclusion; an easy check 
shows that (\ref{orthogonality}) is satisfied, hence our hypothesis (\ref{given}) holds true in this case.\\

\noindent
{\bf Example 3.}
We take $M$ to be a symplectic manifold with moment map $J:M\rightarrow \mathfrak{g}^{*}$; $\rho$ will be the coadjoint representation, and we can
use the Killing form to identify $\mathfrak{g}$ with $\mathfrak{g}^{*}$. 
Again, $\Phi_m$ is restricted to be a linear combination of the $G$-invariants 
in ${\rm Sym}^{2m}(\mathfrak{g})$ and $||v||^{2m}$ is one of them. 
In particular, if $G={\rm SO}(3)$, then 
any invariant must be a scalar multiple of $||v||^{2m}$ and our hypothesis 
(\ref{given}) must hold. A particular case of this situation that 
illustrates the usefulness of our formula (\ref{thesis}) will be 
discussed in the next section.\\

\section{Interacting vortices on a 2-sphere} \label{Vortices}
In this section, we describe an application of Theorem~\ref{thethm} to a
natural setting where our hypothesis holds true.
Let $N$ be a positive integer. We consider the measure space $M$ to be
the moduli space $\mathcal{M}_{N}$ of $N$-vortices on a 2-sphere of radius
$R>\sqrt{N}$~\cite{ManSut}. This is just $\mathbb{CP}^{N}$ as a complex
manifold, but equipped with a K\"ahler structure $\omega_{L^2}$ induced from a
gauge-theoretic version of the $L^2$ norm on the space of fields. The
associated K\"ahler metric $g_{r\bar{s}}$ encodes information about the physics
of vortices at low energies; for example, its geodesic flow gives a good
approximation to the slow dynamics of the abelian Higgs model close to critical
coupling~\cite{StuAH}. For $N>1$, this K\"ahler structure is distinct from the
Fubini--Study structure $\omega_{\rm FS}$ on $\mathbb{CP}^{N}$, although it has
been argued that~\cite{BatMan}
\[
\omega_{L^2} = 2\pi (R^2-N) \omega_{\rm FS} + {o}(R^2 - N)
\]
as $R^2 \searrow {N}$.

There is a local description of $\omega_{L^2}$, which we briefly recall
here~\cite{SamVS, RomGVB}. This uses the fact that 
$\mathbb{CP}^{N}\cong{\rm Sym}^{N}(S^2):=(S^2)^{N}/\mathfrak{S}_N$. We denote
by $\Delta\subset (S^2)^N$ the set of fixed points of elements of
$\mathfrak{S}_N$, and let $z$ be a complex stereographic co\"ordinate on an
open set $U\subset S^{2}$. Then in terms of the natural co\"ordinates
$(z_1,\ldots, z_N)$ for $U^N\subset (S^2)^N$,
\begin{equation}\label{Samols}
\omega_{L^2}\equiv
\frac{i}{2}\sum_{r,s=1}^{N}g_{r\bar{s}}
dz_r \wedge d\bar{z}_s
=i\sum_{r,s=1}^{N}\left( 
\frac{R^2 \delta_{rs}}{(1+|z_r|^2)^2} + \frac{\partial{b}_r}{\partial\bar{z}_s}
\right) dz_r \wedge d\bar{z}_s.
\end{equation}
The functions $b_r(z_1,\ldots,z_N)$ are defined on ${U}^{N}-\Delta$ in terms of
a solution to an elliptic PDE reminiscent of the Liouville equation. They
satisfy
\[
b_r(\ldots,z_r,\ldots, z_s,\ldots)=b_{s}(\ldots,z_s,\ldots,z_r,\ldots),
\]
therefore local quantities like $g_{r\bar{s}}$ in (\ref{Samols}) descend to the
moduli space. Although the $b_r$ are not known explicitly, some statements
about them (and the metric) can be made using the symmetry of the problem. For
example, the fact that ${\rm SO}(3)$ acts on $S^2$ by isometries implies the
relations~\cite{RomQCS}
\[
\sum_{r=1}^{N}(z_r b_r - \bar{z}_r \bar{b}_r)=0 \quad \mbox{and} \quad
\sum_{r=1}^{N}(2z_r + {z_r}^{\!2} \, b_r + \bar{b}_r)=0,
\]
which in turn can be used to show that (\ref{Samols}) is preserved by the
induced action of ${\rm SO}(3)$ on $\mathcal{M}_{N}$. Thus the Liouville
measure associated to $\omega_{L^2}$ is ${\rm SO}(3)$-invariant. In addition,
there exists a moment map 
${J}=(J_1,J_2,J_3):\mathcal{M}_N\rightarrow\mathfrak{so}(3)^{*}$, 
for which our hypothesis (\ref{given}) is obviously satisfied. Its components
can be calculated as~\cite{RomQCS}
\begin{eqnarray*}
J_1 &\!\!\!=&\!\!\! 2 \pi \sum_{r=1}^{N}\left(R^2 \frac{z_r +
\bar{z}_r}{1+|z_r|^2}
+ {\textstyle \frac{1}{2}} (b_r+\bar{b}_r) \right), \\
J_2 &\!\!\!=& \!\!\!-2 \pi i \sum_{r=1}^{N}\left(R^2 \frac{z_r -
\bar{z}_r}{1+|z_r|^2}
- {\textstyle \frac{1}{2}} (b_r-\bar{b}_r) \right), \\
J_3 &\!\!\!=\!\!\!& 2 \pi \sum_{r=1}^{N}\left( R^2 \frac{1-|z_r|^2}{1+|z_r|^2}
- 
(z_r b_r +1) \right).
\end{eqnarray*}

The metric on $\mathcal{M}_N$ has been used to study the statistical mechanics
of a gas of vortices in the abelian Higgs model at critical coupling, both in
the noninteracting case where the net forces experienced by the vortices are
zero~\cite{ManSMV} and in presence of a background potential~\cite{RomGVB}. A
more physically interesting situation would be the case where inter-vortex
interactions are introduced. At the level of the dynamics on the moduli space,
these would be described by a potential $V:\mathcal{M}_N \rightarrow
\mathbb{R}$ invariant under the action of ${\rm SO}(3)$. The simplest
nontrivial potential with this property is just (a multiple of) the square of
the moment map,
\[
V=\mu^2 ||{J}||^{2} = \mu^2 ({J_1}^2 + {J_2}^2 + {J_3}^2),
\]
where $\mu^2$ is just a coupling constant. In the corresponding model, the
partition function for the gas of vortices is given by (see \cite{RomGVB} for
details)
\begin{eqnarray*}
Z& =&\frac{1}{(2 \pi \hbar)^{2N}}\int_{T^{\ast}\mathcal{M}_N}
\exp \left( -({\frac{1}{2\pi}\sum_{r,s}g^{r \bar{s}} w_r\bar{w}_s + \mu^2
||{J}||^2 })/T \right)
\;\frac{\omega_{\rm can}^{2N}}{(2N)!} \\
&=&\left( \frac{T}{2\hbar^2}\right)^{N}\int_{\mathcal{M}_N}
\exp \left( -\frac{\mu^2 }{T}\,
||{J}||^2 \right) \;
\frac{\omega_{L^2}^{N}}{N!}.
\end{eqnarray*}
Here, $w_r=\pi \sum_{s=1}^{N}g_{r\bar{s}} \dot{\bar{z}}_s$ denote canonical
momenta to the moduli $z_r$, 
\[
\omega_{\rm can}=\frac{1}{2}\sum_{r=1}^{N}(dz_r \wedge d\bar{w}_r + 
d\bar{z}_{r} \wedge dw_r)
\]
is the canonical symplectic form on the phase space $T^{*}\mathcal{M}_N$, 
$2 \pi \hbar$ is Planck's constant and $T$ denotes the absolute temperature. It
turns out that Theorem~\ref{thethm} can be used to find a closed expression for
this partition function. In fact, after expanding the exponential, we can
organise the remaining integral as a sum of integrals of powers of the angular
momentum along a fixed axis by making use of (\ref{formulae}) with $n=3$:
\[
Z = \left( \frac{T}{2\hbar^2}\right)^{N}\sum_{m=0}^{\infty}\frac{(-1)^m
\mu^{2m}}{m!\; T^m} (2m+1) \int_{\mathcal{M}_N} {J_{3}}^{2m}.
\]
Now we apply the following result of~\cite{RomGVB}, which is obtained as an
application of the Duistermaat--Heckman formula to the circle action generated
by $J_3$ on the symplectic manifold $(\mathcal{M}_{N},\omega_{L^2})$:
\[
\int_{\mathcal{M}_N} {J_{3}}^{2m} =
\sum_{j=0}^{N}\frac{(-1)^{N-j}(2m)!}{j!(N-j)!(N+2m)!}
\left( 2\pi(R^2 -N) (2j-N)\right)^{N+2m}.
\]
The partition function can then be expressed as
\begin{eqnarray*}
Z&\!\!\!=\!\!\!& \left( \frac{T}{2\hbar^2}\right)^{N}
\sum_{j=0}^{N}\sum_{m=0}^{\infty}\frac{(-1)^{j+m}(2m+1)!}{j!(N-j)!m!(N+2m)!} 
{\textstyle
\left(\frac{\mu^2}{T}\right)^m 
\left( 4 \pi (R^2-N) \left(\frac{N}{2}-j\right)\right)^{N+2m}}\\
&\!\!\!=\!\!\!&\frac{1}{N!}\left( \frac{\tilde{A}T}{2\hbar^2}\right)^N
\sum_{j=0}^{N}
\frac{(-1)^j}{j!(N-j)!} 
{\textstyle
\left(\frac{N}{2}-j \right)^N \!\!\!
\left.\right._{2}\!
F_2 (1,\frac{3}{2};\frac{N+1}{2},
\frac{N+2}{2}; 
-\frac{(\mu \tilde{A})^2} {T} (\frac{N}{2}-j)^2)}
\end{eqnarray*}
in terms of the generalised hypergeometric function 
$\!\left.\right._{2}\!F_2$ \cite{SlapFq}, and where we introduced the area
available for $N$ vortices on the sphere
\[
\tilde{A}:=4 \pi(R^2-N).
\]
It is still a challenging problem to understand the physics determined by $Z$
--- in particular, it would be interesting to obtain an equation of state 
for this system (in some approximation) and analyse whether it allows for 
phase transitions.\\

\begin{small}
\bibliographystyle{numsty}
\bibliography{biblio}
\end{small}

\end{document}